\documentclass[11pt]{article}

\usepackage{graphicx,booktabs}
\usepackage{amsmath}
\usepackage{natbib}
\usepackage{times}
\usepackage{bm}
\usepackage[normalem]{ulem}
\usepackage{color}
\usepackage{xr}
\newtheorem{theorem}{Theorem}
\newtheorem{lemma}{Lemma}
\newcommand{\trans}{^{ \mathrm{\scriptscriptstyle T} }}
\newcommand{\Var}{\text{var}}
\newcommand{\Cov}{\text{cov}}

\newcommand{\effrange}{\textsc{er}}

\begin{document}

\title{The Role of the Range Parameter for Estimation and Prediction in Geostatistics}
\author{C. G. KAUFMAN \and B. A. SHABY}

\maketitle

\begin{abstract}
Two canonical problems in geostatistics are estimating the parameters in a
specified family of stochastic process models and predicting the process at new
locations. A number of asymptotic results addressing these problems over a fixed
spatial domain indicate that, for a Gaussian process with Mat\'ern covariance
function, one can fix the range parameter controlling the rate of decay of the
process and obtain results that are asymptotically equivalent to the case that
the range parameter is known. In this paper we show that the same asymptotic
results can be obtained by jointly estimating both the range and the variance of
the process using maximum likelihood or maximum tapered likelihood. Moreover, we
show that intuition and approximations derived from asymptotic arguments
using a fixed range parameter can be problematic when applied to finite samples,
even for moderate to large sample sizes. In contrast, we show via simulation
that performance on a variety of metrics is improved and asymptotic
approximations are applicable for smaller sample sizes when the range and
variance parameters are jointly estimated. These effects are particularly
apparent when the process is mean square differentiable or the effective range
of spatial correlation is small.

\end{abstract}


\section{Introduction}

The analysis of point-referenced spatial data, often referred to as
geostatistics, relies almost exclusively on a single construct: the stationary
Gaussian process with a parametric mean and covariance function. Exceptions may
be found, some of them notable, but in almost all elaborate hierarchical or
nonstationary models in the literature, one can find structures built from
stationary Gaussian processes.

Given the prominent role of the stationary Gaussian process, it is perhaps
surprising that the theoretical properties of inference under this model remain
incompletely understood. Consider a canonical problem in geostatistics, that of
predicting the value of a spatial process with unknown model parameters at
locations not contained in the dataset. \citet{Stein2010} gives a succinct
overview of asymptotic issues for both parameter estimation and prediction.

\citet[e.g.][]{stein-1999a} makes a compelling case
for using the Mat\'ern covariance model for the Gaussian process $\{Z(s), s \in D \subseteq \Re^d\}$, with
\begin{equation}
  \label{eqn:matern}
  \Cov\{Z(s_i), Z(s_j)\} =
  \sigma^2 K(s_i - s_j; \rho, \nu) = 
  \frac{\sigma^2 (\|s_i - s_j \|/\rho)^\nu}
       {\Gamma(\nu)2^{\nu-1}}
    \mathcal{K}_\nu(\|s_i - s_j \|/\rho),
\end{equation}
where $\sigma^2, \rho, \nu > 0$, and $\mathcal{K}_\nu$ is the modified Bessel
function of the second kind of order $\nu$ \citep[][Section
9$\cdot$6]{Abromowitz1967}. The range parameter $\rho$ controls the rate of
decay with distance, with larger values of $\rho$ corresponding to more highly
correlated observations. This model is particularly attractive because of its
flexibility in representing the smoothness of the Gaussian process, with any
degree of mean square differentiability being possible, according to the value
of $\nu$ \citep{stein-1999a}.

\citet{Zhang2004} provides influential results concerning the consistency of
parameter estimates for the Mat\'ern model under infill, or fixed-domain
asymptotics. Infill asymptotics requires that the sampling domain be fixed as
the number, and hence density, of observations increases to infinity. These
results on consistency follow from a more fundamental result in
\citet{Zhang2004} concerning equivalence, or mutual absolute continuity, of
Gaussian measures on bounded domains.

A theoretical treatment of spatial prediction and corresponding standard error
estimation has been developed in a series of works by
\citet{Stein1988,Stein1990,Stein1993,stein-1999a}. These works provide
conditions under which predictions using a mis-specified covariance function are
asymptotically efficient and associated standard error estimates converge almost
surely to their true values under infill asymptotics. One such condition is that
the mis-specified covariance be chosen so that the resulting Gaussian measure
and the true one are equivalent, providing a link to the results in
\citet{Zhang2004}. However, as we will discuss in Section
\ref{sec:theory-prediction}, the nature of that link has sometimes been
misinterpreted. 

Like most of the aforementioned works, we focus on the isotropic $d$-dimensional
Mat\'ern covariance model \eqref{eqn:matern}. We devote particular attention to
the range parameter $\rho$. One may detect in the recent applied
literature a growing tendency to regard $\rho$ as secondarily influential. For
example, \citet{zhang-2010a} find that fixing $\rho$ at arbitrary large values
has little impact on predictive performance, and \citet{Gneiting2010} argue that
specifying a single $\rho$ for all variables in a multivariate model is not
restrictive. \citet[e.g.][]{sahu-2007a} choose from a small number of fixed
values of $\rho$ based on performance on hold-out data, while
\citet{anderes-2012a} produce spatial predictions without ever estimating $\rho$
from the data. These authors borrow intuition from the asymptotic results of
\citet{Stein1988}, \citet{Zhang2004}, and others, results that fix $\rho$ at an
arbitrary value, often for mathematical tractability. Each of these results
presents some particular variation of the conclusion that fixing $\rho$ at an
incorrect value is asymptotically just as good as using the true value. However,
as we will show, this intuition cannot necessarily be transferred so readily to
the finite sample case. Here, we focus on joint estimation of
$\sigma^2$ and $\rho$, and we prove general results for joint estimation using
maximum likelihood or maximum tapered likelihood under the Mat\'ern model. We
demonstrate via simulation that inference based on these new asymptotic
results is superior on a variety of metrics.

\section{Asymptotic Theory for Estimation and Prediction}
\label{sec:theory-main}

\subsection{Preliminaries}
\label{sec:theory-prelim}

We begin with some notation and assumptions that will be used in all our results unless specifically stated otherwise.  Let $Z = \{Z(s), s \in D \subset \Re^d\}$ be a stochastic process on a bounded domain $D$, with $d = 1, 2,$ or $3$.  Let $G(0, \sigma^2 K_\theta)$ denote the mean zero stationary Gaussian measure for $Z$ with marginal variance $\sigma^2 > 0$ and correlation function $K_\theta$, depending on parameters $\theta \in \Theta \subseteq \Re^p$.  For a particular sampling design $S_n = \{s_1,\ldots,s_n \in D\}$ of distinct locations, we observe $Z_n = \{Z(s_1), \ldots Z(s_n)\}\trans$.  Our tasks are to use $Z_n$ to estimate $\sigma^2$ and $\theta$ and to predict $Z(s_0)$ for some location $s_0 \in D$, not in $S_n$.  Our results concern the behaviour of these estimators and predictors under infill asymptotics.  

We use $G(0, \sigma^2 K_{\rho, \nu})$ to denote a mean zero Gaussian
measure with the Mat\'ern covariance function. We also assume that the
smoothness parameter $\nu$ is known. Our focus is on the role played by the
range parameter $\rho$ in this model, namely to show that several important
results that have been provable only in the case of fixing $\rho$ at an
arbitrary value can be extended to the case that $\rho$ is estimated.

The reason that it is justifiable to fix $\rho$, at least in an asymptotic sense, follows from a property of the Mat\'ern model shown by \citet{Zhang2004}.  This result indicates that when the dimension $d \leq 3$, two Gaussian measures with the same $\nu$ but different values of $\rho$ can in fact be equivalent.  Specifically, Theorem 2 of 
\citet{Zhang2004} states that for fixed $\nu > 0$, $G(0, \sigma^2_0 K_{\rho_0, \nu})$ and $G(0, \sigma^2_1 K_{\rho_1, \nu})$ are equivalent on bounded domains if and only if $\sigma^2_0/\rho_0^{2\nu} = \sigma^2_1/\rho_1^{2\nu}$.

The parameter $c = \sigma^2 / \rho^{2\nu}$ is what \citet{stein-1999a} calls a microergodic parameter. \citet[][page 175]{stein-1999a}  suggests re-parametrizing into microergodic and non-microergodic components of the parameter vector, which we here define as $c$ and $\rho$, respectively.  He conjectures that if all model parameters are estimated by maximum likelihood, the asymptotic behaviour of the maximum likelihood estimator for the microergodic parameter is the same as if the non-ergodic component were known.  In the next section, we outline existing results that concern the asymptotic behaviour for the maximum likelihood estimator for $c$ when $\rho$ is fixed at an arbitrary value, and we extend them to the case that $\rho$ is estimated, showing that Stein's conjecture is true for the Mat\'ern model.

\subsection{Estimation of Covariance Parameters}
\label{sec:theory-estimation}

Theorem 2 of \citet{Zhang2004} has an immediate and important corollary for estimation, namely that 
there do not exist consistent estimators of $\sigma^2$ or $\rho$ based on a sequence of observation vectors taken at an increasing sequence of subsets of a bounded domain.  However, 
it is important to note that this corollary does not imply that the data contain no information about $\sigma^2$ and $\rho$ individually.  Indeed, in simulation studies we observe that sampling distributions for the maximum likelihood estimators can in many cases be quite concentrated about the true values, even as we know these distributions will not become ever more concentrated as $n$ increases \citep{Zhang2004,Kaufman2006}.  Some intuition behind this can be given by appealing to another asymptotic framework, that of increasing the domain of observations.  \citet{Mardia1984} give regularity conditions under which the maximum likelihood estimators for all model parameters are consistent and asymptotically normally distributed, and these conditions may be shown to hold under an increasing domain framework.  Any finite set of observation locations could conceivably be a member in a sequence under either the fixed-domain or increasing-domain asymptotic framework. \citet{Zhang2005} note that the increasing domain framework can be mimicked by fixing the domain but decreasing the range parameter.  Therefore, it is not surprising that when the true range parameter is small relative to the sampling domain, it can be well estimated from data.

The likelihood function for $\sigma^2$ and $\rho$ under the Mat\'ern model with fixed $\nu>0$ based on observations $Z_n$ is
\begin{equation}
\label{eqn:likelihood}
\mathcal{L}_n(\sigma^2, \rho) = (2\pi\sigma^2)^{-n/2} |\Gamma_n(\rho)|^{-1/2} 
\exp\left\{-\frac{1}{2\sigma^2} Z_n\trans \Gamma_n(\rho)^{-1} Z_n\right\},
\end{equation}
where $\Gamma_n(\rho)$ is the matrix with entries $K(s_i-s_j; \rho, \nu)$ $(i,j =1,\ldots,n)$ for $K$ defined as in (\ref{eqn:matern}).  We consider two types of estimators obtained by maximizing (\ref{eqn:likelihood}).  The first fixes $\hat\rho_n = \rho_1$ for all $n$ and maximizes $\mathcal{L}_n(\sigma^2, \rho_1)$.  The second maximizes (\ref{eqn:likelihood}) over both $\sigma^2$ and $\rho$.  In either case, the estimator of $\sigma^2$ may be written as a function of the corresponding estimator of $\rho$.  That is, we may write 
$
\hat\sigma^2_n(\hat\rho_n) = \arg\max_{\sigma^2} \mathcal{L}_n(\sigma^2, \hat\rho_n) = Z_n\trans \Gamma_n(\hat\rho_n)^{-1} Z_n / n,
$
where $\hat\rho_n$ is either $\rho_1$ or the value the maximizes the profile likelihood for $\rho$, when a unique maximizer exists.  In most cases the latter estimator is not available in closed form and must be found numerically.  We may likewise express the corresponding estimators of $c = \sigma^2/\rho^{2\nu}$ as a function of $\hat\rho_n$, namely 
\begin{equation}
\label{eqn:chatdefinition}
\hat{c}_n(\hat\rho_n) = \hat\sigma^2_n(\hat\rho_n)/\hat\rho_n^{2\nu} = Z_n\trans \Gamma_n(\hat\rho_n)^{-1} Z_n / (n\hat\rho_n^{2\nu}).
\end{equation}

We state for reference the following result defining the asymptotic behaviour of $\hat{c}_n(\rho_1)$ for an arbitrary fixed value $\rho_1>0$.  This result combines Theorem 3 of \citet{Zhang2004} and Theorem 3 of \citet{Wang2011}.

\begin{theorem}
\label{thm:fixrhoasymptotics}
  Let $S_n$ be an increasing sequence of subsets of $D$.  Then as $n \rightarrow \infty$, 
  \begin{enumerate}
  \item $\hat{c}_n(\rho_1) \rightarrow \sigma^2_0/\rho_0^{2\nu}$ almost surely, and
  \item $n^{1/2} \{\hat{c}_n(\rho_1) - \sigma^2_0/\rho_0^{2\nu}\} \rightarrow N\{0, 2 (\sigma^2_0/\rho_0^{2\nu})^2\}$ in distribution
  \end{enumerate}
  under $G(0, \sigma^2_0 K_{\rho_0, \nu}).$
\end{theorem}


A key contribution of the current paper is to show that Theorem \ref{thm:fixrhoasymptotics} can be used as a stepping stone to proving that the maximum likelihood estimator $\hat{c}_n(\hat{\rho}_n)$ has exactly the same asymptotic behaviour as does $\hat{c}_n(\rho_1)$ for any $\rho_1$, including the true value $\rho_0$.  We make use of the following lemma, which shows that $\hat{c}_n(\hat\rho_n)$ is monotone when viewed as a function of $\hat\rho_n$. 

\begin{lemma}
  \label{lem:monotonicity}
  Let $S_n = \{s_1, \ldots, s_n \in D \subseteq \Re^d\}$ denote any set of observation locations in any dimension.  Fix $\nu > 0$ and define $\Gamma_n(\rho)$ to be the matrix with entries $K(s_i - s_j; \rho, \nu)$ as in (\ref{eqn:matern}).  Define $\hat{c}_n(\rho) =
  Z_n\trans\Gamma_n(\rho)^{-1}Z_n/(n\rho^{2\nu})$.  Then for any $0 < \rho_1 < \rho_2$, $\hat{c}_n(\rho_2) \leq \hat{c}_n(\rho_1)$ for any vector $Z_n$.
    \end{lemma}

{\em Proof.

  Let $0 < \rho_1 < \rho_2$.  The difference
$$
 \hat{c}_n(\rho_1) - \hat{c}_n(\rho_2) 
=
    Z_n\trans\{\rho_1^{-2\nu}\Gamma_n(\rho_1)^{-1} -
                \rho_2^{-2\nu}\Gamma(\rho_2)^{-1}\}
    Z _n/n
$$
  is non-negative for any $Z_n$ if the matrix 
  $A = \rho_1^{-2\nu}\Gamma_n(\rho_1)^{-1} -\rho_2^{-2\nu}\Gamma_n(\rho_2)^{-1}$ 
  is positive semi-definite.  By Corollary 7$\cdot$7$\cdot$4(a) of \citet[page 473]{horn-1985a},
  $A$ is positive semi-definite if and only if the matrix 
  $B = \rho_2^{2\nu}\Gamma_n(\rho_2)  - \rho_1^{2\nu}\Gamma_n(\rho_1)$ 
  is positive semi-definite.  The entries of $B$ may be expressed in terms of a function $K_B: \Re^d  \rightarrow \Re$, with
  $$
  B_{ij} = K_B(s_i - s_j) = \rho_2^{2\nu} K(\|s_i - s_j\|;
  \rho_2, \nu) - \rho_1^{2\nu} K(\|s_i - s_j\|;
  \rho_1, \nu),
  $$
  and $B$ is positive semi-definite if $K_B$
  is a positive definite function.  Define
    \begin{align}
    f_B(\omega) & = \label{eqn:fourier-ordering}
    \frac{1}{(2\pi)^d} 
    \int_{\Re^d}\!\text{e}^{-i\omega\trans x}
                         K_B(x)\,\mbox{d}x 
    \nonumber \\
    & = 
    \frac{1}{(2\pi)^d} \bigg\{
      \rho_2^{2\nu}
      \int_{\Re^d}\!\text{e}^{-i\omega\trans x}
                           K(x; \rho_2, \nu)\,\mbox{d}x -
      \rho_1^{2\nu}
      \int_{\Re^d}\!\text{e}^{-i\omega\trans x}
                           K(x; \rho_1, \nu)\,\mbox{d}x
    \bigg\}.
  \end{align}
  Both integral terms in (\ref{eqn:fourier-ordering}) are finite, with
  $$
  \frac{1}{(2\pi)^d} \int_{\Re^d} e^{-i \omega\trans x} K(x; \rho, \nu) dx =
  \frac{\Gamma(\nu+d/2)}{\pi^{d/2}\Gamma(\nu)}\rho^{-2\nu}
  (\rho^{-2}+\|\omega\|^2)^{-(\nu+d/2)},$$ 
  the spectral density of the Mat\'ern correlation function.  Therefore,
$$
    f_B(\omega) = \frac{\Gamma(\nu+d/2)}{2^d\pi^{3d/2}\Gamma(\nu)}
    \Big\{(\rho_2^{-2} + \|\omega\|^2)^{-(\nu+d/2)} - 
         (\rho_1^{-2} + \|\omega\|^2)^{-(\nu+d/2)}\Big\}.
$$         
  To show $K_B$ is positive definite it suffices to show $f_B(\omega)$ is positive for all $\omega$.  This is clear because $0 < \rho_1 < \rho_2$.   Therefore $\hat{c}_n(\rho_2) \leq \hat{c}_n(\rho_1)$ for any vector $Z_n$.}

We can now make use of Theorem \ref{thm:fixrhoasymptotics} in proving a more general result for the maximum likelihood estimator when the parameter space for $\rho$ is taken to be a bounded interval. This condition was also used by \citet{Ying1991}, who proved Theorem \ref{thm:estrhoasymptotics} for the special case that $D$ is the unit interval and $\nu = 1/2$.  These bounds are not restrictive in practice, as the interval may be taken to be arbitrarily large.

\begin{theorem}
\label{thm:estrhoasymptotics}
  Let $S_n$ be an increasing sequence of subsets of $D$.
  Suppose $(\sigma_0^2, \rho_0)\trans \in (0, \infty) \times [\rho_L, \rho_U]$, for
  any $0 < \rho_L < \rho_U < \infty$. Let $(\hat{\sigma}^2_n, \hat{\rho}_n)\trans$ maximize (\ref{eqn:likelihood}) over $(0, \infty)
  \times [\rho_L, \rho_U]$. Then 
  \begin{enumerate} 
  \item $\hat{\sigma}^2_n/ \hat{\rho}_n^{2\nu}
  \rightarrow \sigma^2_0/\rho_0^{2\nu}$ almost surely, and
  \item $n^{1/2} (\hat\sigma^2_n/\hat\rho_n^{2\nu} - \sigma^2_0/\rho_0^{2\nu}) \rightarrow N\{0, 2 (\sigma^2_0/\rho_0^{2\nu})^2\}$ in distribution
  \end{enumerate}
  under $G(0, \sigma^2_0 K_{\rho_0,\nu})$.
\end{theorem}

{\em Proof. 

  By assumption, $\rho_L \leq \hat{\rho}_n \leq \rho_U$ for every $n$. Define two sequences, $\hat{c}_n(\rho_L)$ and $\hat{c}_n(\rho_U)$, according to (\ref{eqn:chatdefinition}).  By Lemma \ref{lem:monotonicity}, $\hat{c}_n(\rho_L) \leq \hat{c}_n(\hat\rho_n) = \hat\sigma^2_n/\hat{\rho}_n^{2\nu} \leq \hat{c}_n(\rho_U)$ for all $n$ with probability one.  Combining this with Theorem \ref{thm:fixrhoasymptotics} applied to $\hat{c}_n(\rho_L)$ and $\hat{c}_n(\rho_U)$ implies the result.}
 
Theorem \ref{thm:estrhoasymptotics} is useful because it applies to
the procedure that is most often adopted in practice, of allowing the range
parameter to be estimated from data over some bounded domain. In fact,
the method of proof in Theorem \ref{thm:estrhoasymptotics} implies that these
asymptotic results hold for any bounded sequence $\hat\rho_n$, provided that
$\hat\sigma^2_n$ is defined as in (\ref{eqn:chatdefinition}). This would
include, for example, estimating $\rho$ using the variogram and plugging it into
(\ref{eqn:chatdefinition}), but not joint estimation of $\rho$ and $\sigma^2$
using the variogram. In practice, the bounds for numerical optimization of
$\rho$ can be chosen to be arbitrarily wide, subject to numerical stability.

A similar method of proof can be used to show consistency and asymptotic
normality of the maximum tapered likelihood estimator proposed by
\citet{Kaufman2008}. The online supplement contains analogues of Lemma
\ref{lem:monotonicity} and Theorem \ref{thm:estrhoasymptotics} for this estimator. 

Arguments following from \citet{Zhang2004} would suggest that the range
parameter may be fixed in practice. However, as we shall show in Section
\ref{sec:simulation-main}, the estimator $\hat{c}_n(\rho_1)$ can often display
sizeable bias, making the approximation in Theorem \ref{thm:fixrhoasymptotics}
quite inaccurate. Confidence intervals constructed using Theorem
\ref{thm:fixrhoasymptotics} can, due to this bias, have empirical coverage
probabilities very near to zero in some cases. In contrast, we will show that
confidence intervals for $c$ constructed using Theorem
\ref{thm:estrhoasymptotics} have close to nominal coverage even for moderate
sample sizes.

\subsection{Prediction at New Locations}
\label{sec:theory-prediction}

We now consider the problem of predicting the value of the process at a new location $s_0$ not in the set of observation locations $S_n$.  \citet{Stein1988, Stein1990, Stein1993, stein-1999a} has considered this problem when an incorrect model is used.  
Predictors under the wrong model can be consistent under relatively weak conditions.  Our focus is therefore on two other desirable properties, asymptotic efficiency and asymptotically correct estimation of prediction variance.  In a seminal paper, \citet{Stein1988} showed that both of these properties hold when the model used is equivalent to the true measure. In the case of the Mat\'ern covariance, Theorem 2 of \citet{Zhang2004} indicates that this holds for a model with the correct $\nu$ and microergodic parameter $\sigma^2/\rho^{2\nu}$.  This has led to statements in the literature to the effect that the parameter c = $\sigma^2/\rho^{2\nu}$ can be consistently estimated, and this is what matters for prediction.  While this statement contains an element of truth, we will argue in this section that it can also be somewhat misleading, both in an asymptotic sense, as well as in guiding choices for applications.

Under the mean zero Gaussian process model with Mat\'ern covariance function and known $\nu>0$, define
\begin{equation}
  \label{eqn:krig-mean}
\hat{Z}_n(\rho) = \gamma_n(\rho)\trans \Gamma_n(\rho)^{-1} Z_n,
\end{equation}
where $\gamma_n(\rho) = \{K(s_0 - s_i; \rho,\nu)\}_i$ and $\Gamma_n(\rho) =
\{K(s_i - s_j; \rho,\nu)\}$ ($i,j = 1,\ldots, n$). The predictor
$\hat{Z}_n(\rho)$ is the best linear unbiased predictor for $Z_0 = Z(s_0)$ under
a presumed model $G(0, \sigma^2 K_{\rho, \nu})$ for any value of $\sigma^2$.
This predictor does not depend on $\sigma^2$, only $\rho$ and $\nu$. Therefore,
any intuition that one can fix $\rho = \rho_1$, and that plug-in predictions
will improve with $n$ due in any way to convergence of $\hat{c}_n(\rho_1)$ with
$n$, is clearly a misunderstanding of asymptotic results. Equivalence, although
sufficient for asymptotic efficiency, is not necessary. The way in which $c$ is
relevant for prediction concerns estimates of the mean squared error of the
predictor. Under model $G(0, \sigma^2_0 K_{\rho_0, \nu})$, this is
\begin{align}
  \label{eqn:krig-variance}
  \Var_{\sigma^2_0, \rho_0}\{\hat{Z}_n(\rho)-Z_0\} = \sigma^2_0 
    & \{ 1 - 2 \gamma_n (\rho)\trans\Gamma_n(\rho)^{-1}\gamma_n(\rho_0) +  \\
    \nonumber 
    & \gamma_n(\rho)\trans \Gamma_n(\rho)^{-1} \Gamma_n(\rho_0) 
      \Gamma_n(\rho)^{-1} \gamma_n(\rho)\},
\end{align}
where $\gamma_n(\rho_0)$ and $\Gamma_n(\rho_0)$ are defined analogously to their
counterparts using $\rho$. In the case that $\rho = \rho_0$, this expression
simplifies to
\begin{equation}
  \label{eqn:krig-variance-naive}
  \Var_{\sigma^2_0, \rho_0}\{\hat{Z}_n(\rho_0)-Z_0\} = 
  \sigma^2_0 \{1 - \gamma_n(\rho_0)\trans\Gamma_n(\rho_0)^{-1}\gamma_n(\rho_0)\}.
\end{equation}
In practice, it is common to estimate the model parameters and then plug them into (\ref{eqn:krig-mean}) and (\ref{eqn:krig-variance-naive}), treating them as known. The asymptotic properties of this procedure, so-called plug-in prediction, are quite difficult to obtain.  Instead, most theoretical development has been under a framework in which plug-in parameters are fixed, rather than being estimated from observations at an increasing sequence of locations. We will review these results and indicate how they may be extended to include estimation of the variance parameter $\sigma^2$ with a fixed value of $\rho$, making precise the sense in which the statement regarding $c$ at the beginning of this section should be interpreted.

The following result is an application of Theorems 1 and 2 of \citet{Stein1993}. 

\begin{theorem}
  \label{thm:stein-1993}
  Suppose $G(0, \sigma^2_0 K_{\rho_0, \nu})$ and $G(0, \sigma^2_1 K_{\rho_1, \nu})$ are two Gaussian process measures on $D$ with the same value of $\nu  > 0$.
  \begin{enumerate}
  \item  As $n \rightarrow \infty$,
  $$
  \frac{\Var_{\sigma^2_0, \rho_0}\big\{\hat{Z}_n(\rho_1) - Z_0\big\}}
       {\Var_{\sigma^2_0, \rho_0}\big\{\hat{Z}_n(\rho_0) - Z_0\big\}} 
  \rightarrow 1.
  $$
  \item Furthermore, if $\sigma^2_0/\rho_0^{2\nu} = \sigma^2_1/\rho_1^{2\nu}$, then as $n \rightarrow \infty$,
   \begin{equation}
   \label{eqn:varratio1}
  \frac{\Var_{\sigma^2_1, \rho_1}\big\{\hat{Z}_n(\rho_1) - Z_0\big\}}
       {\Var_{\sigma^2_0, \rho_0}\big\{\hat{Z}_n(\rho_1) - Z_0\big\}} 
  \rightarrow 1.
  \end{equation}
  \end{enumerate}
\end{theorem}

{\em Proof.

Let $f_0$ be the spectral density corresponding to $\sigma^2_0 K_{\rho_0, \nu}$ and $f_1$ be the spectral density corresponding to $\sigma^2_1 K_{\rho_1, \nu}$.  The result follows from noting that the function $f_0(\omega) \| \omega \| ^{2\nu + d}$ is bounded away from zero and infinity as $\| \omega \| \rightarrow \infty$ and that
$$
\lim_{\| \omega \| \rightarrow \infty} \frac{f_1(\omega)}{f_0(\omega)} = \frac{\sigma^2_1 / \rho_1^{2\nu}}{\sigma^2_0 / \rho_0^{2\nu}}.
$$
These two conditions satisfy those needed for Theorems 1 and 2 of \citet{Stein1993}.}

The implication of part 1 of Theorem \ref{thm:stein-1993} is that if the correct value of $\nu$ is used, any value of $\rho$ will give asymptotic efficiency.  The condition $\sigma^2_0/\rho_0^{2\nu} = \sigma^2_1/\rho_1^{2\nu}$ is not necessary for asymptotic efficiency, but it does provide asymptotically correct estimates of mean squared prediction error.  The numerator in (\ref{eqn:varratio1}) is the naive mean squared error for $\hat{Z}_n(\sigma^2_1, \rho_1)$, assuming model $G(0, \sigma^2_1 K_{\rho_1, \nu})$, whereas the denominator is the true mean squared error for $\hat{Z}_n(\sigma^2_1, \rho_1)$, under model $G(0, \sigma^2_0 K_{\rho_0, \nu})$.  We now show the same convergence happens if $\rho$ is fixed at $\rho_1$ but $\sigma^2$ is estimated via maximum likelihood.  This is an extension of part 2 of Theorem \ref{thm:stein-1993}.  Part 1 needs no extension, since the form of the predictor itself does not depend on $\sigma^2$.

\begin{theorem}
  \label{thm:consistency-pred-var}
  Suppose $G(0, \sigma^2_0 K_{\rho_0, \nu})$ is a Gaussian process measure on $D$.  Fix $\rho_1 > 0$.  For a sequence of observations $Z_n$ on an increasing sequence of subsets $S_n$ of $D$, define $\hat\sigma^2_n = Z_n\trans \Gamma_n(\rho_1)^{-1} Z_n / n$. Then as $n \rightarrow \infty$,
\begin{equation}
\label{eqn:varratio2}
  \frac{\Var_{\hat\sigma^2_n, \rho_1}\big\{\hat{Z}_n(\rho_1) - Z_0\big\}}
       {\Var_{\sigma^2_0, \rho_0}\big\{\hat{Z}_n(\rho_1) - Z_0\big\}} \rightarrow 1
\end{equation}
almost surely under $G(0, \sigma^2_0 K_{\rho_0, \nu})$.
\end{theorem}

{\em Proof.

Define $\sigma^2_1 = \sigma^2_0(\rho_1/\rho_0)^{2\nu}$.  Then write
$$ 
\frac{\Var_{\hat\sigma^2_n, \rho_1}\big\{\hat{Z}_n(\rho_1) - Z_0\big\}}
     {\Var_{\sigma^2_0, \rho_0}\big\{\hat{Z}_n(\rho_1) - Z_0\big\}} = 
\frac{\Var_{\hat\sigma^2_n, \rho_1}\big\{\hat{Z}_n(\rho_1) - Z_0\big\}}
     {\Var_{\sigma^2_1, \rho_1}\big\{\hat{Z}_n(\rho_1) - Z_0\big\}}
\frac{\Var_{\sigma^2_1, \rho_1}\big\{\hat{Z}_n(\rho_1) - Z_0\big\}}
     {\Var_{\sigma^2_0, \rho_0}\big\{\hat{Z}_n(\rho_1) - Z_0\big\}}.
$$
By Theorem \ref{thm:stein-1993}, 
$\Var_{\sigma^2_1, \rho_1}\{\hat{Z}_n(\rho_1) - Z_0\}/\Var_{\sigma^2_0, \rho_0}\{\hat{Z}_n(\rho_1) - Z_0\} \rightarrow 1.$  
So we need only show that 
$\Var_{\hat\sigma^2, \rho_1}\{\hat{Z}_n(\rho_1) - Z_0\}/\Var_{\sigma^2_1, \rho_1}\{\hat{Z}_n(\rho_1) - Z_0\} \rightarrow 1$ almost surely under $G(0, \sigma^2_0 K_{\rho_0, \nu})$.  By (\ref{eqn:krig-variance-naive}),
$\Var_{\hat\sigma^2_n, \rho_1}\{\hat{Z}_n(\rho_1) - Z_0\}/\Var_{\sigma^2_1, \rho_1}\{\hat{Z}_n(\rho_1) - Z_0\} = \hat\sigma^2_n/\sigma^2_1$.  Under $G(0, \sigma^2_1 K_{\rho_1, \nu}),$ $\hat\sigma^2_n$ is equal in distribution to $\sigma^2_1/n$ times a $\chi^2$ random variable with $n$ degrees of freedom and hence converges almost surely to $\sigma^2_1$ as $n \rightarrow \infty$.  Because $\sigma^2_0/\rho_0^{2\nu} = \sigma^2_1/\rho_1^{2\nu},$ Theorem 2 of \citet{Zhang2004} gives that $G(0, \sigma^2_0 K_{\rho_0, \nu})$ and $G(0, \sigma^2_1 K_{\rho_1, \nu})$ are equivalent, so that $\hat\sigma^2_n \rightarrow \sigma^2_1$ almost surely under $G(0, \sigma^2_0 K_{\rho_0, \nu})$ as well.}

We conjecture that the asymptotic behaviour in Theorem \ref{thm:stein-1993}$\cdot$1 and Theorem \ref{thm:consistency-pred-var} still holds if $\rho_1$ is replaced by $\hat\rho_n$, the maximum likelihood estimator, although proving such a result has thus far been intractable for cases of practical interest \citep{Putter2001}.

\section{Simulation Study}
\label{sec:simulation-main}

\subsection{Setup}
\label{sec:simulation-setup}Œ

Fixing the range parameter is supported by asymptotic results, and it is computationally efficient in practice.  However, it is unclear to what degree asymptotic results are appropriate in guiding our choices for applied problems with finite sample sizes.  To systematically explore this issue, we simulate data under a Gaussian process model for a variety of settings chosen to mimic the range of behaviour we might observe in practice, and we compare the performance of inference procedures that either fix or estimate the range parameter.

\begin{figure}
   \centering
   \includegraphics[height=0.8\textwidth,angle=-90]{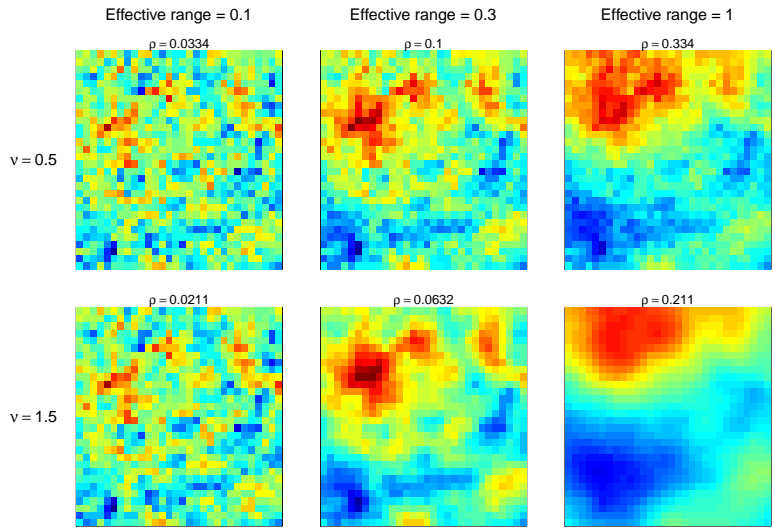} 
   \caption{Simulated random fields on $[0,1]^2$ under parameter settings used in the simulation study.  Fields were simulated using the Cholesky method, using the same set of random deviates for each panel.  The value of the range parameter $\rho$ corresponding to each $\nu$ and effective range combination is also indicated.}
   \label{fig:fieldsillustration}
\end{figure}

We simulate data in the unit square under the mean zero Gaussian process model
with Mat\'ern covariance and smoothness parameter $\nu = 0\cdot 5$ or 1$\cdot$5
and marginal variance $\sigma^2=1$. We also use three effective ranges for the
process, choosing values of $\rho$ such that the correlation decays to
0$\cdot$05 at distances of 0$\cdot$1, 0$\cdot$3, or 1. Figure
\ref{fig:fieldsillustration} illustrates the effect of these parameter settings.
As we shall see, whether a particular sample size is large enough such that
finite sample properties are well approximated by asymptotic results depends on
both the degree of smoothness and the effective range of the process.  

We also vary the sample size in the simulation, taking $n=400, 900$, and $1600$.  To avoid numerical issues that can arise from sampling locations situated too close to each other, sampling locations are constructed using a perturbed grid.  We construct a $67 \times 67$ regular grid with coordinates from 0$\cdot$005 to 0$\cdot$995 in increments of 0$\cdot$015 in each dimension.  To each gridpoint, we add a random perturbation according to a uniform distribution over [-0$\cdot$005, 0$\cdot$005]$^2$.  The resulting set of 4489 locations therefore has the property that all pairs of points have at least 0$\cdot$005 distance from each other.  We then choose random subsets of these locations to be our observation locations, with each sample size containing the points from smaller sample sizes.  In evaluating the predictive properties of models fit using a fixed or estimated range parameter, we consider a $50\times 50$ regular grid of locations over $[0,1]^2$.  

For each parameter setting, we simulate $1000$ datasets corresponding to
realizations of the Gaussian process observed at the union of $n = 1600$
observation and $m=2500$ prediction locations. For each dataset and sample size,
we estimate $\sigma^2$ and $\rho$ by numerically maximizing the profile
likelihood for $\rho$ and plugging the result into the corresponding closed-form
estimator for $\sigma^2$.
We also calculate $\hat\sigma^2_n(\rho_1) = Z_n\trans
\Gamma_n(\rho_1)^{-1}Z_n/n$ for values of $\rho_1$ equal to 0$\cdot$2,
0$\cdot$5, 1, 2, and 5 times the true value of $\rho$. Corresponding to each of
these parameter estimates, we also construct 95\% confidence intervals for $c =
\sigma^2/\rho^{2\nu}$ using the normal approximation provided by Theorem
\ref{thm:fixrhoasymptotics} when $\rho$ is fixed and Theorem
\ref{thm:estrhoasymptotics} when $\rho$ is estimated. Finally, we construct
kriging predictors and estimated standard errors for each of the $m=2500$
prediction locations by plugging in parameter estimates into
(\ref{eqn:krig-mean}) and (\ref{eqn:krig-variance-naive}).

Optimization was carried out using the R function optim with the L-BFGS-B
option, which we restricted to the interval $\rho \in [\varepsilon, 15\rho_0]$,
where $\varepsilon$ is defined by machine precision, about $10^{-16}$ on our
machine. Neither endpoint was ever returned.


In the following sections we discuss the results for estimation and prediction.  Many of the results show a similar pattern, which can be summarized as follows.  Not surprisingly, the performance of the maximum likelihood estimator, maximizing over both $\sigma^2$ and $\rho$, is generally very good, especially by $n=1600$.  Procedures using a fixed $\rho$ are almost always worse, although there are certain settings under which the differences are minimal.  These tend to be for $\nu=1/2$ (corresponding to processes that are not mean-square differentiable) and a large effective range.  In these cases, particularly when $\rho$ is fixed at something larger than its true value, the estimators and predictors can still perform well.  This agrees with some examples in the literature, for which $\nu=1/2$ and large effective ranges were used \citep{zhang-2010a, Wang2011}.  When the process is smooth ($\nu$=1$\cdot$5) and/or the true range of spatial correlation is small, estimation and prediction is markedly improved by estimating $\rho$ via maximum likelihood.

\subsection{Parameter estimation}
\label{sec:simulation-estimation}

\begin{figure}
   \centering
   \includegraphics[height=0.9\textwidth,angle=-90]{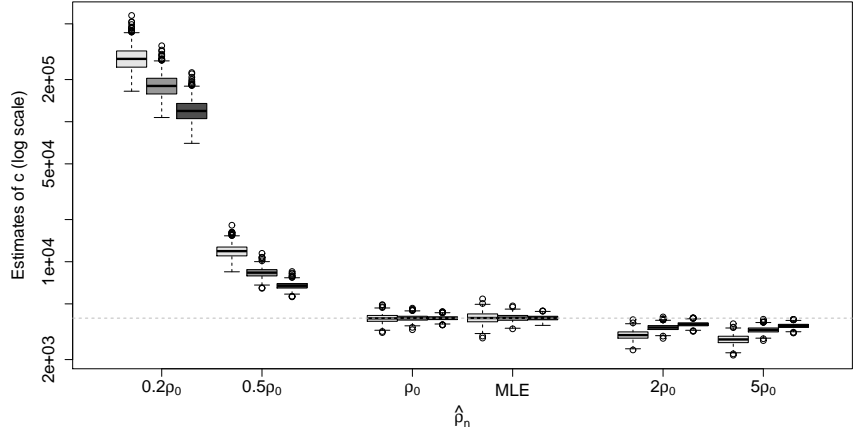} 
   \caption{Sampling distributions for $\hat{c}_n$ when $\nu$=1$\cdot$5 and the
           effective range is 0$\cdot$3. The range parameter is either fixed at the
            true value ($\rho_0$), estimated via maximum likelihood (MLE), or fixed at a multiple of the
            truth (0$\cdot2\rho_0$, $\ldots$ , $5\rho_0$). The four boxplots in each
            group correspond to sample sizes of $n=400$, 900, and 1600,
            reading from left to right.}
   \label{fig:cboxplots}
\end{figure}

Given the asymptotic results in \citet{Zhang2004} and \citet{Wang2011} for
$\hat{c}_n(\rho_1)$ for fixed $\rho_1$, it is tempting to adopt the intuition
that this estimator can adapt to incorrectly specified values of $\rho$. While
this is true asymptotically, we observe in our simulation results that this
adaptation in many cases requires a very large value of $n$; sampling
distributions can be highly biased and can move very slowly toward the truth as
$n$ increases. Figure \ref{fig:cboxplots} illustrates these effects for a subset
of our simulation results, namely when $\nu$=1$\cdot$5 and the effective range
is 0$\cdot$3. Sampling distributions for $\hat{c}_n(\rho_1)$ are noticeably
biased. As we expect from Theorem \ref{thm:fixrhoasymptotics}, these biases
decrease with $n$, although even when $n=1600$ the true value of $c$ lies far in
the tail of the sampling distribution. In contrast, the sampling distributions
for the maximum likelihood estimator $\hat{c}_n(\hat\rho_n)$ have negligible
bias. Indeed, they behave very similarly to those for the estimator of $c$ that
fixes $\rho$ at the truth. Similar effects can be seen for other values of $\nu$
and effective range.  See Table \ref{tbl:cbiasreltruth} in the supplement for the
relative bias of different estimators of $c$. 

\begin{table}[ht] \small
\caption{Empirical coverage rates of confidence intervals for $c = \sigma^2/\rho^{2\nu}$, expressed as percentages.  Intervals are constructed using either the maximum likelihood estimator (MLE) or estimates of $c$ that fix $\hat\rho_n$ at a multiple of the true value of $\rho$, rounded to the nearest 1\%.\label{tbl:ccoverage}}
\begin{center}
\begin{tabular}{crrrrrrr}
&& \multicolumn{3}{c}{$\nu$=0$\cdot$5} & \multicolumn{3}{c}{$\nu$=1$\cdot$5} \\
$\hat\rho_n$ &$n$ & \effrange=0$\cdot$1 & \effrange =0$\cdot$3 & \effrange=1 & \effrange=0$\cdot$1 & \effrange=0$\cdot$3 & \effrange=1\\
  \addlinespace[0.5em]
  MLE 
   & 400 & 81 & 92 & 94 & 64 & 87 & 94 \\ 
   & 900 & 89 & 94 & 94 & 74 & 91 & 94 \\ 
   & 1600 & 90 & 94 & 94 & 81 & 92 & 95 \\ 
  \addlinespace[0.5em]
  0$\cdot$2 $\rho$ 
   & 400 & 0 & 0 & 0 & 0 & 0 & 0 \\ 
   & 900 & 0 & 0 & 1 & 0 & 0 & 0 \\ 
   & 1600 & 0 & 0 & 2 & 0 & 0 & 0 \\ 
  \addlinespace[0.5em]
  0$\cdot$5 $\rho$ 
   & 400 & 0 & 4 & 88 & 0 & 0 & 4 \\ 
   & 900 & 0 & 7 & 90 & 0 & 0 & 9 \\ 
   & 1600 & 0 & 13 & 92 & 0 & 0 & 18 \\
  \addlinespace[0.5em] 
  2 $\rho$ 
   & 400 & 3 & 75 & 93 & 0 & 1 & 83 \\ 
   & 900 & 3 & 82 & 93 & 0 & 9 & 89 \\ 
   & 1600 & 5 & 84 & 94 & 0 & 17 & 93 \\ 
  \addlinespace[0.5em]
  5 $\rho$ 
   & 400 & 0 & 63 & 92 & 0 & 0 & 77 \\ 
   & 900 & 0 & 75 & 93 & 0 & 2 & 86 \\ 
   & 1600 & 0 & 79 & 93 & 0 & 5 & 90 \\
\end{tabular}
\end{center}
\end{table}

If Theorem \ref{thm:fixrhoasymptotics} is used to construct confidence intervals
and $n$ is in fact not large enough for the normal approximation to be
appropriate, the coverage of such intervals can be disastrously low. Table
\ref{tbl:ccoverage} shows empirical coverage rates for confidence intervals
constructed as $\hat{c}_n(\hat\rho_n) \pm 1$$\cdot$$96 \,
\{2\hat{c}_n(\hat\rho_n)^2/n\}^{1/2}$ for $\hat\rho_n$ equal to the maximum
likelihood estimator or a fixed $\rho_1$. Theorems \ref{thm:fixrhoasymptotics}
and \ref{thm:estrhoasymptotics} imply that these intervals are asymptotically
valid 95\% confidence intervals for $c$. Not surprisingly, however, given the
large biases observed when $\rho$ is fixed, the differences in the empirical
coverage rates between fixed and estimated $\rho$ are striking, even when $n$ is
large. In many cases the coverage for intervals constructed using
$\hat{c}_n(\rho_1)$ was 0\%, to within Monte Carlo sampling error. Coverage is
best when $\nu$ is small and effective range is large. For fixed $\rho_1$, it
also appears to be better to choose $\rho_1 > \rho_0$ than $\rho_1 < \rho_0$.

\subsection{Prediction}
\label{sec:simulation-prediction}

The mean squared error of predictor $\hat{Z}_n(\rho_1)$ may be calculated in closed form using (\ref{eqn:krig-variance}).  When the plug-in predictor $\hat{Z}_n(\hat\rho_n)$ is used, we need to integrate over the sampling distribution for $\hat\rho_n$, which we approximate by averaging over the simulation results from Section \ref{sec:simulation-estimation}. For both fixed and estimated $\rho$, we calculate the average mean squared prediction error, averaging over the $m=2500$ prediction points in the domain.  Because the prediction problem varies in difficulty according to $\nu$, effective range, and sample size $n$, we report the percent increase in mean squared prediction error relative to the optimal mean squared prediction error using the true value of $\rho$, which is calculated from (\ref{eqn:krig-variance-naive}).  These results are shown in Table \ref{tbl:incmse}.

\begin{table}[ht] \small
\caption{Percent increase in mean squared prediction error relative to the optimal mean squared prediction error using the true value of $\rho$, rounded to the nearest 0$\cdot$1 percent.\label{tbl:incmse}}
\begin{center}
\begin{tabular}{crrrrrrr}
&& \multicolumn{3}{c}{$\nu$=0$\cdot$5} & \multicolumn{3}{c}{$\nu$=1$\cdot$5} \\
$\hat\rho_n$ &$n$ & \effrange=0$\cdot$1 & \effrange=0$\cdot$3 & \effrange=1 & \effrange=0$\cdot$1 & \effrange=0$\cdot$3 & \effrange=1\\
  \addlinespace[0.5em]
  MLE 
   & 400 & 0$\cdot$2 & 0$\cdot$1 & 0$\cdot$0 & 0$\cdot$2 & 0$\cdot$1 & 0$\cdot$1 \\ 
   & 900 & 0$\cdot$1 & 0$\cdot$0 & 0$\cdot$0 & 0$\cdot$1 & 0$\cdot$0 & 0$\cdot$0 \\ 
   & 1600 & 0$\cdot$0 & 0$\cdot$0 & 0$\cdot$0 & 0$\cdot$0 & 0$\cdot$0 &  0$\cdot$0\\ 
   \addlinespace[0.5em]
   0$\cdot$2 $\rho$ 
   & 400 & 36$\cdot$6 & 60$\cdot$4 & 6$\cdot$5 & 103$\cdot$1 & 487$\cdot$0 & 165$\cdot$5 \\ 
   & 900 & 56$\cdot$4 & 37$\cdot$5 & 2$\cdot$3 & 218$\cdot$2 & 474$\cdot$1 & 83$\cdot$8 \\ 
   & 1600 & 66$\cdot$2 & 19$\cdot$2 & 0$\cdot$9 & 351$\cdot$4 & 321$\cdot$5 & 41$\cdot$5 \\ 
   \addlinespace[0.5em]
   0$\cdot$5 $\rho$ 
   & 400 & 8$\cdot$7 & 2$\cdot$8 & 0$\cdot$2 & 26$\cdot$9 & 20$\cdot$0 & 2$\cdot$9 \\ 
   & 900 & 7$\cdot$9 & 1$\cdot$1 & 0$\cdot$1 & 32$\cdot$9 & 10$\cdot$2 & 1$\cdot$3 \\ 
   & 1600 & 5$\cdot$5 & 0$\cdot$4 & 0$\cdot$0 & 29$\cdot$2 & 4$\cdot$7 & 0$\cdot$7 \\ 
   \addlinespace[0.5em]
   2 $\rho$ 
   & 400 & 2$\cdot$8 & 0$\cdot$3 & 0$\cdot$0 & 12$\cdot$0 & 2$\cdot$1 & 0$\cdot$3 \\ 
   & 900 & 1$\cdot$3 & 0$\cdot$1 & 0$\cdot$0 & 6$\cdot$8 & 1$\cdot$0 & 0$\cdot$1 \\ 
   & 1600 & 0$\cdot$6 & 0$\cdot$0 & 0$\cdot$0 & 3$\cdot$4 & 0$\cdot$4 & 0$\cdot$1 \\ 
   \addlinespace[0.5em]
   5 $\rho$ 
   & 400 & 5$\cdot$6 & 0$\cdot$6 & 0$\cdot$1 & 27$\cdot$2 & 4$\cdot$2 & 0$\cdot$6 \\ 
   & 900 & 2$\cdot$4 & 0$\cdot$2 & 0$\cdot$0 & 13$\cdot$7 & 1$\cdot$9 & 0$\cdot$2 \\ 
   & 1600 & 1$\cdot$1 & 0$\cdot$1 & 0$\cdot$0 & 6$\cdot$6 & 0$\cdot$9 & 0$\cdot$1 \\
\end{tabular}
\end{center}
\end{table}

It is clear from Table \ref{tbl:incmse} that plug-in prediction using the maximum likelihood estimator $\hat\rho_n$ performs quite well relative to predicting with the true value of $\rho$.  For $n=900$ and 1600, the increase in mean squared error is less than 0$\cdot$1 percent in all cases.  It is also clear that there are cases in which it makes little difference if $\rho$ is fixed at an incorrect value.  This is true when the effective range is large and $\rho_1$ is fixed at something larger than the true value.  However, there are also cases in which fixing $\rho$ can lead to quite a large loss of efficiency.  These effects are magnified when we move from $\nu=$0$\cdot$5 to $\nu$=1$\cdot$5, suggesting that a misspecified value of $\rho$ is more problematic for smoother processes.  This aligns with some earlier cases in the literature in which predictions with a fixed $\rho$ were still quite accurate.  For example, \citet{zhang-2010a} examined precipitation data using a predictive process model \citep{Banerjee2008} and concluded that a variety of prediction metrics did not change when $\rho$ was fixed at values larger than the maximum likelihood estimator.  However, the underlying covariance model for the predictive process was Mat\'ern with $\nu$=0$\cdot$5, corresponding to a process that is not mean square differentiable.

In a similar pattern to what we observe for mean squared error in Table \ref{tbl:incmse}, using the maximum likelihood estimator produces intervals with the nominal rate in nearly all cases, and the estimators fixing $\rho$ at something larger than the true value achieve this rate for $n = 900$ and 1600 when the effective range is large.  However, the intervals tend to be too conservative when the effective range is large and $\rho_1$ is too small, and they tend to be not conservative enough when the effective range is small and $\rho_1$ is too big.  See the supplement for full results.

\section{Discussion}
\label{sec:discussion}

We have made a number of simplifying assumptions. Considering the ways in which
these assumptions may be relaxed provides a rich set of questions for future
research. For example, our results concern only mean zero Gaussian processes,
which is equivalent to assuming that the mean of the process is known. Results
on equivalence of mean zero Gaussian measures such as Theorem 2 of
\citet{Zhang2004} can be used in proving equivalence of Gaussian process
measures with different means \citep[][Chapter 4, Corollary 5]{stein-1999a}.
However, the primary difficulty is in extending estimation results.
\citet{Zhang2004} indicates that his method of proof is not easily extended to
the case of an unknown mean term. Asymptotic results for the case $\nu=1/2$ and
$d=1$ are given in Theorem 3 of \citet{Ying1991}, and it seems plausible that
similar results might hold for $d=2$ and 3. With an unknown mean, it might be
preferable to use restricted maximum likelihood \citet{stein-1999a}, for which
improved infill asymptotic results should also be sought.

We have also not considered what happens when the observations are not of the
process $Z$ itself, but of $Z$ observed with measurement error. Again, results
for equivalence and prediction can be extended in a relatively straightforward
way. We expect something like Theorem \ref{thm:estrhoasymptotics} should hold
for the case that $Z$ is observed with measurement error. However, in a
restricted version of this problem, the introduction of the error term reduces
the rate of convergence of the maximum likelihood estimator for $c$ from the
usual order $n^{-1/2}$ to order $n^{-1/4}$ \citep{Chen2000}.  

Perhaps the most important restriction, both here and in previous work, is that
the smoothness parameter $\nu$ is assumed to be known. Estimating $\nu$ provides
desirable flexibility, as this parameter controls the mean square
differentiability of the process. However, we know of no results concerning the
maximum likelihood estimator in this case. \citet[][Section
6$\cdot$7]{stein-1999a} examines a periodic version of the Mat\'ern model and
argues that $\hat\sigma^2_n$ and $\hat\nu_n$ should have a joint asymptotic
normal distribution, but it is an open question whether a similar result holds
for non-periodic fields. 
\section*{Acknowledgements}
Cari Kaufman's portion of this work was supported by the Center for Science of Information (CSoI), an NSF Science and Technology Center, under grant agreement CCF-0939370.  Benjamin Shaby's portion of this work was supported by NSF Grant DMS-06-35449 to the Statistical and Applied Mathematical Sciences Institute. The authors thank the editor, associate editor, and two anonymous reviewers for their useful suggestions.


\bibliography{RangeParameter}

\begin{thebibliography}{21}
\expandafter\ifx\csname natexlab\endcsname\relax\def\natexlab#1{#1}\fi

\bibitem[{Abromowitz \& Stegun(1967)}]{Abromowitz1967}
\textsc{Abromowitz, M.} \& \textsc{Stegun, I.}, eds. (1967).
\newblock \textit{Handbook of Mathematical Functions}.
\newblock U.S. Government Printing Office.

\bibitem[{Anderes et~al.(2012)Anderes, Huser, Nychka \& Coram}]{anderes-2012a}
\textsc{Anderes, E.}, \textsc{Huser, R.}, \textsc{Nychka, D.} \& \textsc{Coram,
  M.} (2012).
\newblock Nonstationary positive definite tapering on the plane.
\newblock \textit{J. Comput. Graph. Statist.} To appear.

\bibitem[{Banerjee et~al.(2008)Banerjee, Gelfand, Finley \&
  Sang}]{Banerjee2008}
\textsc{Banerjee, S.}, \textsc{Gelfand, A.}, \textsc{Finley, A.} \&
  \textsc{Sang, H.} (2008).
\newblock Gaussian predictive process models for large spatial data sets.
\newblock \textit{J. Roy. Statist. Soc. Ser. C} \textbf{70}, 825--848.

\bibitem[{Chen et~al.(2000)Chen, Simpson \& Ying}]{Chen2000}
\textsc{Chen, H.}, \textsc{Simpson, D.} \& \textsc{Ying, Z.} (2000).
\newblock Infill asymptotics for a stochastic process model with measurement
  error.
\newblock \textit{Statist. Sinica} \textbf{10}, 141--156.

\bibitem[{Gneiting et~al.(2010)Gneiting, Kleiber \& Schlather}]{Gneiting2010}
\textsc{Gneiting, T.}, \textsc{Kleiber, W.} \& \textsc{Schlather, M.} (2010).
\newblock Mat{\'e}rn cross-covariance functions for multivariate random fields.
\newblock \textit{J. Am. Statist. Assoc.} \textbf{105}, 1167--1177.

\bibitem[{Horn \& Johnson(1985)}]{horn-1985a}
\textsc{Horn, R.~A.} \& \textsc{Johnson, C.~R.} (1985).
\newblock \textit{Matrix analysis}.
\newblock Cambridge: Cambridge University Press.

\bibitem[{Kaufman(2006)}]{Kaufman2006}
\textsc{Kaufman, C.} (2006).
\newblock \textit{Covariance tapering for likelihood-based estimation in large
  spatial data sets}.
\newblock Ph.D. thesis, Carnegie Mellon University.

\bibitem[{Kaufman et~al.(2008)Kaufman, Schervish \& Nychka}]{Kaufman2008}
\textsc{Kaufman, C.}, \textsc{Schervish, M.} \& \textsc{Nychka, D.} (2008).
\newblock Covariance tapering for likelihood-based estimation in large spatial
  data sets.
\newblock \textit{J. Am. Statist. Assoc.} \textbf{103}, 1545--1555.

\bibitem[{Mardia \& Marshall(1984)}]{Mardia1984}
\textsc{Mardia, K.} \& \textsc{Marshall, R.} (1984).
\newblock Maximum likelihood estimation of models for residual covariance in
  spatial regression.
\newblock \textit{Biometrika} \textbf{71}, 135--146.

\bibitem[{Putter \& Young(2001)}]{Putter2001}
\textsc{Putter, H.} \& \textsc{Young, G.} (2001).
\newblock On the effect of covariance function estimation on the accuracy of
  kriging predictors.
\newblock \textit{Bernoulli} \textbf{7}, 421--438.

\bibitem[{Sahu et~al.(2007)Sahu, Gelfand \& Holland}]{sahu-2007a}
\textsc{Sahu, S.~K.}, \textsc{Gelfand, A.~E.} \& \textsc{Holland, D.~M.}
  (2007).
\newblock High-resolution space-time ozone modeling for assessing trends.
\newblock \textit{J. Am. Statist. Assoc.} \textbf{102}, 1221--1234.

\bibitem[{Stein(1988)}]{Stein1988}
\textsc{Stein, M.} (1988).
\newblock Asymptotically efficient prediction of a random field with a
  misspecified covariance function.
\newblock \textit{Ann. Statist.} \textbf{16}, 53--63.

\bibitem[{Stein(1990)}]{Stein1990}
\textsc{Stein, M.} (1990).
\newblock Uniform asymptotic optimality of linear predictions of a random field
  using an incorrect second-order structure.
\newblock \textit{Ann. Statist.} \textbf{18}, 850--872.

\bibitem[{Stein(1993)}]{Stein1993}
\textsc{Stein, M.} (1993).
\newblock A simple condition for asymptotic optimality of linear predictions of
  random fields.
\newblock \textit{Statist. Probab. Lett.} \textbf{17}, 399--404.

\bibitem[{Stein(2010)}]{Stein2010}
\textsc{Stein, M.} (2010).
\newblock Asymptotics for spatial processes.
\newblock In \textit{Handbook of Spatial Statistics}, A.~Gelfand, P.~Diggle,
  M.~Fuentes \& P.~Guttorp, eds. Boca Raton, FL: CRC Press, pp. 79--88.

\bibitem[{Stein(1999)}]{stein-1999a}
\textsc{Stein, M.~L.} (1999).
\newblock \textit{Interpolation of Spatial Data: Some theory for kriging}.
\newblock Springer Series in Statistics. New York: Springer-Verlag.

\bibitem[{Wang \& Loh(2011)}]{Wang2011}
\textsc{Wang, D.} \& \textsc{Loh, W.} (2011).
\newblock On fixed-domain asymptotics and covariance tapering in {G}aussian
  random field models.
\newblock \textit{Electron. J. Stat.} \textbf{5}, 238--269.

\bibitem[{Ying(1991)}]{Ying1991}
\textsc{Ying, Z.} (1991).
\newblock Asymptotic properties of a maximum likelihood estimator with data
  from a gaussian process.
\newblock \textit{J. Multivariate Anal.} \textbf{36}, 280--296.

\bibitem[{Zhang(2004)}]{Zhang2004}
\textsc{Zhang, H.} (2004).
\newblock Inconsistent estimation and asymptotically equal interpolations in
  model-based geostatistics.
\newblock \textit{J. Am. Statist. Assoc.} \textbf{99}, 250--261.

\bibitem[{Zhang \& Wang(2010)}]{zhang-2010a}
\textsc{Zhang, H.} \& \textsc{Wang, Y.} (2010).
\newblock Kriging and cross-validation for massive spatial data.
\newblock \textit{Environmetrics} \textbf{21}, 290--304.

\bibitem[{Zhang \& Zimmerman(2005)}]{Zhang2005}
\textsc{Zhang, H.} \& \textsc{Zimmerman, D.} (2005).
\newblock Towards reconciling two asymptotic frameworks in spatial statistics.
\newblock \textit{Biometrika} \textbf{92}, 921--936.

\end{thebibliography}

\end{document}